\documentclass [a4paper] {article}

\usepackage[utf8]{inputenc}

\usepackage[french]{babel}

\usepackage{amssymb,amsmath, amsfonts}

\usepackage [dvips] {graphicx,color}

\usepackage{amsthm}

\usepackage{hyperref}

\usepackage{MnSymbol}

\newcommand*{\twosquig}{%
\mathrel{\vcenter{\offinterlineskip
\hbox{$\rightsquigarrow$}\vskip-.08ex\hbox{$\rightsquigarrow$}}}}

\newtheorem{prop} {Proposition}

\theoremstyle{definition}
\newtheorem{df}{Définition} 
\newtheorem*{df*}{Définition}

\theoremstyle{remark}
\newtheorem{rmq}{Remarque} 
\newtheorem{example}{Exemple} 
\newtheorem{exm}[example]{Exemple}

\author{Stéphane \textsc{Dugowson}
\footnote{Laboratoire Quartz - Institut Polytechnique Grand Paris - Supméca. Email : s.dugowson@gmail.com}
}

\title {Interaction des dynamiques graphiques ouvertes \\ (définitions)}
\date{7 août 2015 (version 3)} 

\begin{document}

\maketitle

\noindent \textbf{Abstract}. The aim of this paper is to define what we shall call \emph{open graphic dynamics}, their interactions and the dynamics produced by those interactions. It prepares the study of  \emph{open sub-categorical dynamics} and  \emph{open categorical dynamics}.

\mbox{}

Relations. Multiple binary relations. Connectivity structures.

\mbox{}

\noindent\textbf{Résumé}. L'objet du présent article est de poser les définitions relatives à ce que nous appellerons des \emph{dynamiques graphiques ouvertes}, aux interactions entre ces dynamiques et aux nouvelles dynamiques engendrées par ces interactions, et cela dans le but de constituer un socle pour l'étude ultérieure des \emph{dynamiques sous-catégoriques ouvertes} et, parmi elles, des \emph{dynamiques catégoriques ouvertes}, et de leurs interactions.

\mbox{}

\noindent \textit{Mots clés}:  Interaction.  Dynamiques. Dynamiques graphiques. Relations. Relations binaires multiples. Structures connectives. 

\mbox{}


\paragraph{Avertissement.} Par rapport aux versions antérieures, cette troisième version du texte précise notamment la notion de dynamique ouverte graphique engendrée par une famille dynamique graphique en distinguant plusieurs types d'engendrements.

\section*{Introduction}

En définissant les \emph{dynamiques graphiques ouvertes}, leurs interactions et les nouvelles dynamiques engendrées par ces interactions, nous souhaitons avec cet article constituer un socle pour l'étude ultérieure de l'interaction entre \emph{dynamiques catégoriques ouvertes}, étude dont la principale difficulté vient de ce que, dans le cas général, l'interaction de dynamiques catégoriques ouvertes ne produit pas une dynamique qui soit encore catégorique\footnote{Les difficultés rencontrées tiennent essentiellement à la composition des transitions en tant qu'elles expriment des possibilités pour les dynamiques considérées. 
\`{A} ce sujet, voir notre exposé  \cite{Dugowson:20150506}.
}. 

Du reste, s'agissant déjà des seules dynamiques graphiques, le présent article vise seulement à poser les définitions, et pratiquement aucun exemple ne sera donné ici, réservant pour un article ultérieure la présentation d'un certains nombre d'exemples et l'étude de situations intéressantes.

La section \ref{section preliminaires} donne les définitions préalables et les notations  dont nous aurons besoin, portant respectivement sur les \emph{relations multiples}, sur les \emph{relations binaires multiples}, les transitions et les graphes.

La section \ref{section dynamiques graphiques ouvertes} définit les \emph{dynamiques graphiques ouvertes}. La section \ref{section familles dynamiques} définit les notions d'\emph{interaction} et de \emph{famille dynamique}, c'est-à-dire de famille de dynamiques ouvertes synchronisées en interaction. La section \ref{section DyGrOuv engendrees} définit enfin diverses dynamiques graphiques ouvertes \emph{engendrées} par une famille dynamique.

Le texte se termine par une conclusion intitulée \og Et maintenant ?\fg, les références bibliographiques et une table des matières.

\section{Notions préliminaires et notations employées}
\label{section preliminaires}

\subsection{Relations multiples}

Nous rappelons  la définition des relations multiples,  renvoyant à \cite{Dugowson:201505} en particulier 
pour la définition du monoïde commutatif qui y était noté $(\mathcal{R}_{\mathcal{E}}, \bowtie, 1)$ mais que nous noterons dorénavant  
$(\mathcal{R}_{\subset\mathcal{E}}, \otimes, 1)$,
ainsi que pour la définition de la structure connective d'une relation multiple que nous appliquerons aux relations binaires multiples dans la section \ref{subsec relations binaires multiples}).

\subsubsection{Définition des relations multiples}
\label{subsubsec definition des relations multiples}

 Étant donnée  $\mathcal{E}=(E_i)_{i\in I}$ une famille d'ensembles, nous désignerons  le produit  $\prod_{i\in I} E_i$  par $\Pi_I(\mathcal{E})$ ou $\Pi_I{\mathcal{E}}$, voire simplement par $\Pi_I$ s'il n'y a pas d'ambiguïté sur le contexte.
 
 Plus généralement, pour toute partie $J\subset I$, nous noterons 
\[\Pi_J(\mathcal{E})=\prod_{i\in J} E_i,\] et, s'il n'y a pas d'ambiguïté sur le contexte, nous écrirons simplement $\Pi_J$ plutôt que $\Pi_J(\mathcal{E})$ ou $\Pi_J\mathcal{E}$.

\begin{df}\label{df relation multiple} Une relation multiple $R$ est la donnée d'un triplet $R=(I,\mathcal{E},G)$ constitué
\begin{itemize}
\item d'un ensemble $I$, appelé \emph{index}, \emph{multiplicité}, \emph{arité}, ou encore \emph{domaine} de la relation $R$,
\item d'une famille $\mathcal{E}=(E_i)_{i\in I}$ d'ensembles indexée par $I$, famille appelée \emph{contexte de $R$}, et telle que $\Pi_I\mathcal{E}\neq \emptyset$,
\item d'une partie $G\subset \Pi_I\mathcal{E}$, appelée \emph{graphe de $R$}.
\end{itemize}
\end{df}

\paragraph{Notation.} Nous noterons $\mathcal{R}_I$ la classe constituée de toutes les relations de multiplicité $I$.

\begin{exm} En notant $2$ un ensemble ordonné à deux éléments, la classe $\mathcal{R}_2$ s'identifie à celle des relations binaires entre ensembles.
\end{exm}

\begin{df}[Relations multiples \emph{dans} le contexte $\mathcal{E}$]
L'ensemble $I$ et le contexte $\mathcal{E}=(E_i)_{i\in I}$ étant donnés, nous noterons
$\mathcal{R}_{\subset\mathcal{E}}$ l'ensemble des relations multiples \emph{dans} ce contexte $\mathcal{E}$, c'est-à-dire l'ensemble des relations multiples de la forme $(J,\mathcal{F},H)$ avec
 $J\subset I$, 
 $\mathcal{F}=\mathcal{E}_{\vert J}=(E_j)_{j\in J}$ et
$H\subset \Pi_J=\Pi_J(\mathcal{E})=\Pi_J(\mathcal{F})$. 
 \end{df}

Ainsi,  on a
\[\mathcal{R}_{\subset\mathcal{E}} \subset \bigcup_{J\subset I}\mathcal{R}_J.\]

Pour tout $J\subset I$, les $J$-relations dans $\mathcal{E}$ sont ordonnées par l'inclusion de leurs graphes. Étant données deux $J$-relations $R$ et $S$, nous écrirons $R\subset_J S$, ou simplement $R\subset S$, pour exprimer le fait que $G_R\subset G_S$. 

\subsubsection{Exemples : relations nulles et relations triviales}

La $J$-relation minimale dans $\mathcal{E}$, notée $0_J$, est celle de graphe vide :
\[0_J=(J,\mathcal{E},\emptyset),\]
tandis que la $J$-relation maximale dans le même contexte $\mathcal{E}$, notée $1_J$, est celle de graphe total
\[1_J=(J,\mathcal{E},\Pi_J(\mathcal{E})).\]

Les relations de la forme $0_J$ seront dites \emph{nulles}, celles de la forme $1_J$ seront dites \emph{triviales} (ou \emph{totales}).

Si l'ensemble $J$ est fini, ou par exemple si l'axiome du choix s'applique, alors $\Pi_J\neq\emptyset$, de sorte que $0_J\neq 1_J$. Ceci est vrai, bien que ce ne soit pas très intuitif, en particulier pour $J=\emptyset$, auquel cas le graphe de $0_\emptyset$ est vide tandis que celui de $1_\emptyset$ est $\Pi_\emptyset=\{\bullet\}$.
Cette dernière relation  \og pleine... sur aucun ensemble\fg\, est l'élément neutre 
\[\mathbf{1}=1_\emptyset=(\emptyset, \{\bullet\})\]
de la loi de composition binaire commutative $\otimes$ du monoïde $(\mathcal{R}_{\subset\mathcal{E}}, \otimes, 1)$ (voir \cite{Dugowson:201505}, section § \textbf{1.5.1.}).

\subsection{Relations binaires multiples de domaine $I$}\label{subsec relations binaires multiples}

\subsubsection{Définition}
\label{subsubsec df relations binaires multiples}
\begin{df}\label{df relation binaire multiple} Une \emph{relation binaire multiple} $C$ est un quadruplet $(I,\mathcal{A},\mathcal{B},G)$ où
\begin{itemize}
\item $I$ est un ensemble, appelé \emph{multiplicité},  \emph{domaine}, \emph{arité}, \emph{uplicité} ou encore \emph{index}  de $C$, 
\item $\mathcal{A}=(A_i)_{i\in I}$ est une famille d'ensembles indexée par $I$, appelée \emph{contexte d'entrée de $C$},
\item $\mathcal{B}=(B_i)_{i\in I}$ est une famille indexée par $I$ d'ensembles non vides, appelée \emph{contexte de sortie de $C$},
\item $G\subset \Pi_I\mathcal{E}=\prod_{i\in I}E_i$ est le \emph{graphe} de $C$, où $\mathcal{E}=(E_i=A_i\times B_i)_{i\in I}$ est le \emph{contexte de $C$}, les contextes d'entrée et de sortie étant tels que $\Pi_I\mathcal{E}\neq\emptyset$.
\end{itemize}
\end{df}

\begin{rmq}\label{rmq designation contexte} La donnée du triplet $(I,\mathcal{A}, \mathcal{B})$ équivaut à celle du contexte $\mathcal{E}$, aussi le contexte de $C$ pourra-t-il être désigné aussi bien comme étant $\mathcal{E}$, $(\mathcal{A}, \mathcal{B})$, $(I,\mathcal{A}, \mathcal{B})$ ou encore $(I,\mathcal{A}, \mathcal{B},\mathcal{E})$.
\end{rmq}

\paragraph{Notations.} Nous noterons 
\begin{itemize}
\item $\mathcal{BM}$ la classe de toutes les relations binaires multiples,
\item $\mathcal{BM}_I$ la classe des relations binaires multiples de multiplicité $I$,
\item $\mathcal{BM}_{(\mathcal{A}, \mathcal{B})}$, $\mathcal{BM}_{\mathcal{E}}$ ou $\mathcal{BM}_{(I,\mathcal{A}, \mathcal{B},\mathcal{E})}$ la classe des relations binaires multiples de contexte 
$(I,\mathcal{A}, \mathcal{B},\mathcal{E})$.
\end{itemize}

 .

\begin{rmq} L'hypothèse $\Pi_I\mathcal{E}\neq\emptyset$ entraîne également la non vacuité de $\Pi_I\mathcal{A}$ et de $\Pi_I\mathcal{B}$, ce qui entraîne à son tour que, pour tout $i\in I$, $A_i\neq\emptyset$ et $B_i\neq\emptyset$. Puisque nous admettons l'axiome du choix, on pourrait se contenter de ces dernières hypothèses. Du reste, en pratique, les ensembles $A_i$ et $B_i$ auxquels on fera appel seront tels que même si l'on ne supposait pas l'axiome du choix on aurait toujours
\[(\forall i\in I, A_i\neq\emptyset\neq B_i)\Rightarrow \Pi_I\mathcal{A}\neq\emptyset\neq \Pi_I\mathcal{B}.\]
\end{rmq}

\subsubsection{Trois injections canoniques}\label{subsubsec trois injections}

\paragraph{Injection $rd:\mathcal{BM}_I \hookrightarrow \mathcal{R}_{2I}$.}
Pour tout index $I$, nous noterons $2I$ l'ensemble défini par \[2I=I\sqcup I=I\times\{0,1\}. \]

Pour toute relation binaire multiple $C$ de contextes respectifs $\mathcal{A}=(A_i)_{i\in I}$ et $\mathcal{B}=(B_i)_{i\in I}$, posons 
\[\mathcal{D}=(D_k)_{k\in 2I},\]
où, pour tout $i\in I$, on prend $D_{(i,0)}=A_i$  et $D_{(i,1)}=B_i$. 
Appelons en outre $rd_\mathcal{E}$ l'application 
de $\Pi_I\mathcal{E}$ 
dans $\Pi_{2I}\mathcal{D}$ 
définie pour tout $\left((a_i,b_i)\right)_{i\in I}\in\Pi_I\mathcal{E} $ par
\[rd_\mathcal{E}(\left((a_i,b_i)\right)_{i\in I})=(d_k)_{k\in 2I},\]où pour tout $i\in I$ on prend $d_{(i,0)}=a_i$ et $d_{(i,1)}=b_i$. 

On définit alors une injection canonique  
\[rd:\mathcal{BM}_I \hookrightarrow \mathcal{R}_{2I}\]  en posant, pour toute relation binaire multiple $C=(I,\mathcal{A},\mathcal{B},G)\in \mathcal{BM}_I$,
\[rd(C)=(2I,\mathcal{D}, rd_\mathcal{E}(G))\in \mathcal{R}_{2I},\]
où $rd_\mathcal{E}(G)=\{rd_\mathcal{E}(u), u\in G\}$.

\paragraph{Injection $rm:\mathcal{BM}_I\hookrightarrow \mathcal{R}_I$.}

On définit canoniquement  une injection $rm:\mathcal{BM}_I\hookrightarrow \mathcal{R}_I$ de  la classe des relations binaires multiples de domaine $I$ dans celle des relations multiples également de domaine $I$ en posant
\[rm(I,\mathcal{A},\mathcal{B},G)=(I,\mathcal{E},G),\]
où, comme dans la définition \ref{df relation binaire multiple}, on prend $\mathcal{E}=(A_i\times B_i)_{i\in I}$.

\paragraph{Injection $rb:\mathcal{BM}_I\hookrightarrow \mathcal{R}_2$.}
On définit canoniquement  une injection $rb:\mathcal{BM}_I\hookrightarrow \mathcal{R}_2$ de  la classe des relations binaires multiples de domaine $I$ dans celle des relations binaires en définissant pour toute relation binaire multiple $C=(I,\mathcal{A},\mathcal{B},G)$ la relation binaire $rb(C)$ par
\begin{itemize}
\item sa source : $\Pi_I\mathcal{A}$,
\item son but : $\Pi_I\mathcal{B}$,
\item son graphe : $rb_\mathcal{E}(G)$,
\end{itemize}
où $rb_\mathcal{E}(G)=\{rb_\mathcal{E}(u), u\in G\}$ est l'image de l'ensemble $G$ par l'application $rb_\mathcal{E}:\Pi_I\mathcal{E}\rightarrow\Pi_I\mathcal{A}\times \Pi_I\mathcal{B}$ définie pour tout $(a_i,b_i)_{i\in I} \in \Pi_I\mathcal{E}$ par
\[rb_\mathcal{E}((a_i,b_i)_{i\in I})=((a_i)_{i\in I},(b_i)_{i\in I})\in \Pi_I\mathcal{A}\times \Pi_I\mathcal{B}.\]

\begin{rmq}\label{rmq notation binaire rbC a}
Pour expliciter une relation binaire multiple $C$, il est souvent  commode de donner le graphe $rb_\mathcal{E}(G)$ de la relation binaire associé, en donnant pour toute famille $a\in \Pi_I\mathcal{A}$ l'ensemble des familles $b\in \Pi_I\mathcal{B}$ telles que $(a,b)\in rb_\mathcal{E}(G)$, que l'on pourra noter $rb(C)(a)$ :
\[\forall a\in \Pi_I\mathcal{A}, rb(C)(a)=\{b\in \Pi_I\mathcal{B}, (a,b)\in rb_\mathcal{E}(G)\}.\]
Conformément à la remarque \ref{rmq transitions et relations binaires} faite plus loin page \pageref{rmq transitions et relations binaires}, l'\emph{image} de la relation binaire $rb(C)$ est alors
\[\mathrm{Im}(rb(C))
=\bigcup_{a\in\Pi_I\mathcal{A}} rb(C)(a)
=\{b\in \Pi_I\mathcal{B}, \exists a\in\Pi_I\mathcal{A},(a,b)\in rb_\mathcal{E}(G),\}\] autrement dit, en désignant selon un abus d'écriture usuel le graphe d'une relation binaire par cette relation elle-même,
\[\mathrm{Im}(rb(C))
=
\{b\in \Pi_I\mathcal{B}, \exists a\in\Pi_I\mathcal{A},(a,b)\in rb(C)\}.
\]
\end{rmq}

\begin{rmq}\label{rmq notation rbBmoins1} Pour toute relation binaire $B:\mathcal{A}\rightsquigarrow \Lambda$, et pour tout $\mu\in \Lambda$, nous notons 
\[B^{-1}(\mu)=\{\mathfrak{a}\in\mathcal{A},B(\mathfrak{a},\mu)\}.\] En particulier, pour toute relation binaire multiple $R\in \mathcal{BM}_{(\mathcal{S}, \mathcal{L})}$ avec $\mathcal{S}=(\mathcal{S}_{A_i})_{i\in I}$ et $\mathcal{L}= (L_i)_{i\in I}$, et pour tout $\mu\in\Pi_I(L_i)$, l'expression $rb(R)^{-1}(\mu)$ désigne l'ensemble
\[rb(R)^{-1}(\mu)=\{\mathfrak{a}\in\Pi_I(\mathcal{S}_{A_i}), (\mathfrak{a},\mu)\in rb(R)\}.\]
\end{rmq}

\subsubsection{Structure connective d'une relation binaire multiple}

La définition suivante fait appel à la définition 17 page 15 de \cite{Dugowson:201505}, c'est-à-dire à la définition de la structure connective d'une relation multiple de domaine $I$.

\begin{df} Étant donné un ensemble $I$ et une relation binaire multiple $C$ de domaine $I$, nous appellerons \emph{structure connective de $C$}, et nous noterons $\mathcal{K}_C$, la structure connective de la relation multiple $rm(C)$.
\end{df}

Autrement dit, $\mathcal{K}_C$ est l'ensemble des parties de $I$ non scindables pour la relation multiple $rm(C)$, c'est-à-dire l'ensemble des parties $J$ de $I$ qui n'admettent pas de bipartition $J=K\cup L$ telles que $R_{\vert J}=R_{\vert K}\otimes R_{\vert L}$, où $R=rm(C)$ est la relation multiple définie par $C$.

\subsection{Transitions}\label{subsec notations transitions}

Reprenant les notations de notre article \cite{Dugowson:201112} et de l'ouvrage \cite{Dugowson:201203} --- où sont introduites les dynamiques catégoriques --- nous noterons $\mathcal{P}(A)$ ou simplement $\mathcal{P}A$ l'ensemble des parties de l'ensemble $A$. 
En outre, $\mathbf{P}$  désignera la catégorie
\begin{itemize}
\item dont les objets sont les ensembles,
\item telle que, pour tout couple d'ensembles $(A,B)$, les flèches $f$ de $A$ vers $B$, notées $f:A\rightsquigarrow B$ et appelées \emph{transitions} de $A$ vers $B$, sont les applications de $A$ vers $\mathcal{P}B$,
\item et telle que la composée de deux transitions $f:A\rightsquigarrow B$ et $g:B\rightsquigarrow C$, notée $g\odot f$, est donnée pour tout $a\in A$ par 
\[ g\odot f (a)=\bigcup_{b\in f(a)} g(b). \]
\end{itemize}

\begin{rmq}\label{rmq transitions et relations binaires} La catégorie $\mathbf{P}$ est isomorphe à la catégorie dite \og des relations\fg, c'est-à-dire celle dont les objets sont les ensembles et qui a pour flèches les relations binaires entre eux. Elle admet l'ensemble vide pour objet terminal (objet nul), et le produit cartésien y coïncide avec l'union disjointe. Remarquons d'ailleurs que cette identification entre relation binaire de graphe $f\subset A\times B$ et transition $f:A\rightsquigarrow B$ permet de définir l'\emph{image} $\mathrm{Im} f$ d'une telle relation binaire comme l'image de la transition, c'est-à-dire
\[\mathrm{Im} f = \bigcup_{a\in A}f(a)=\{b\in B,\exists a\in A, (a,b)\in f\}.\]
\end{rmq}

\begin{df}[Transitions (quasi-)déterministes]\label{df transitions quasi deterministes}
Une transition $f:A\rightsquigarrow B$ est dite \emph{quasi-déterministe} si pour tout $a\in A$, $f(a)$ est soit vide, soit un singleton. Dans ce cas, on identifiera $f$ à la fonction $f:A \rightarrow B$ dont le domaine de définition est constitué des $a\in A$ tels que $f(a)\neq\emptyset$. En particulier, $f$ est dite \emph{déterministe} si 
\[\forall a\in A, card(f(a))=1,\] et une telle transition s'identifie à une application que l'on notera encore $f:A\rightarrow B$.
\end{df}

\paragraph{Familles de transitions.}
Pour tout ensemble non vide $M$, nous désignons par $\mathbf{P}^{\underrightarrow{M}}$ la catégorie  dont les objets sont les mêmes que ceux de $\mathbf{P}$ --- autrement dit les ensembles --- 
et dont les flèches de $A$ vers $B$ sont les $M$-familles de transitions de $A$ vers $B$, la composition se faisant composante par composante.

Une famille indexée par $M$ de transitions $(f_\mu:A \rightsquigarrow B)_{\mu\in M}$ de même source $A$ et même but $B$ sera également désignée par l'expression
\[
f:A\twosquig_M B,
\] ou plus simplement, s'il n'y a pas d'ambiguïté sur les paramètres, 
\[
f:A\twosquig B.
\]
La composée de deux telles familles $f:A\twosquig B$ et $g:B\twosquig C$ sera donc la famille 
\[g\odot f:A\twosquig C\] telle que, pour tout $\mu \in M$, on ait
\[(g\odot f)_{\mu}=g_\mu\odot f\mu.\]

\subsection{Graphes}\label{subsec notations graphes}

Les \emph{graphes} dont il est question dans cet article sont des graphes orientés. Un tel graphe $\mathbf{G}$ est défini  par la donnée d'une classe $\dot{\mathbf{G}}$ de sommets, d'une classe $\overrightarrow{\mathbf{G}}$ d'arêtes,
 et pour chaque arête $a\in\overrightarrow{\mathbf{G}}$ de la précision de son sommet source $dom(a)$ et de son sommet but $cod(a)$. 
 Un \emph{morphisme de graphes} $\alpha : \mathbf{F}\to \mathbf{G}$  
 --- que nous appellerons également un \og foncteur de graphes\fg ---
  est la donnée $\alpha=(\dot{\alpha},\overrightarrow{\alpha})$ de deux applications, 
 l'une $\dot{\alpha}:\dot{\mathbf{F}}\to \dot{\mathbf{G}}$ associant des sommets aux sommets, 
 l'autre $\overrightarrow{\alpha}:\overrightarrow{\mathbf{F}}\to \overrightarrow{\mathbf{G}}$ des arêtes aux arêtes, de façon cohérente avec les applications source et but au sens où pour toute arête $\forall a\in\overrightarrow{\mathbf{F}}$, on a
\[
dom(\overrightarrow{\alpha}(a))=\dot{\alpha}(dom(a))
\quad\mathrm{et}\quad 
cod(\overrightarrow{\alpha}(a))=\dot{\alpha}(cod(a)).
 \]

Pour toute catégorie $\mathbf{C}$, nous noterons $Gr(\mathbf{C})$ le graphe associé, dont les sommets sont les objets de $\mathbf{C}$ et dont les arêtes sont les flèches de $\mathbf{C}$.

Par exemple, $Gr(\mathbf{P}^{\underrightarrow{{M}}})$ désigne le graphe dont les sommets sont les ensembles et dont chaque arête est une famille indexée par $M$ de transitions entre les sommets de l'arête.

\section{Dynamiques graphiques ouvertes}
\label{section dynamiques graphiques ouvertes}

\subsection{Dynamiques graphiques}

Comme évoqué au début de la section \ref{subsec notations transitions}, les définitions suivantes généralisent --- et affaiblissent, en les appuyant sur des graphes plutôt que sur des catégories --- celles relatives aux dynamiques catégoriques, aux dynamorphismes, aux horloges, etc. qui ont été données dans \cite{Dugowson:201112} et \cite{Dugowson:201203}.

\subsubsection{Dynamiques, états, déterminisme}

\begin{df}\label{df dynamiques graphiques}[Dynamiques graphiques] Étant donné $\mathbf{G}$ un graphe, une \emph{dynamique graphique} $\alpha$ sur $\mathbf{G}$ est un morphisme de graphes $\alpha:\mathbf{G} \rightarrow Gr(\mathbf{P})$ tel que les images de deux sommets distincts soient deux ensembles disjoints
\[\forall (S,T)\in \dot{\mathbf{G}}^2, S\neq T \Rightarrow \alpha(S)\cap \alpha(T)=\emptyset.\]
Le graphe $\mathbf{G}$ sera appelé le \emph{moteur} de $\alpha$.
\end{df}

Une dynamique $\alpha$ sur $\mathbf{G}$ est ainsi la donnée de deux applications, que l'on pourrait noter respectivement $\dot{\alpha}$ et $\overrightarrow{\alpha}$ mais qu'en pratique on notera le plus souvent toutes deux $\alpha$, qui
\begin{itemize}
\item à tout sommet $S$ de $\mathbf{G}$ associe un ensemble $\dot{\alpha}(S)$ que l'on notera également $S^\alpha$ et dont les éléments seront appelés les \emph{états} de type $S$,
\item à toute arête $a\in\overrightarrow{\mathbf{G}}$ de domaine $S$ et de codomaine $T$ associe une transition $\overrightarrow{\alpha}(a)=a^\alpha:S^\alpha\rightsquigarrow T^\alpha$, autrement dit une application $a^\alpha:S^\alpha\rightarrow\mathcal{P}(T^\alpha)$.
\end{itemize}

\paragraph{Ensemble des états.} Noté $st(\alpha)$, \emph{l'ensemble des états} de la dynamique $\alpha$ est défini par l'union (disjointe)
\[
st(\alpha)=\bigcup_{S\in \dot{\mathbf{G}}}S^\alpha.
\]

Pour tout état $s\in st(\alpha)$, nous noterons $typ(s)$ son type, autrement l'unique sommet $S\in\dot{\mathbf{G}}$ tel que $s\in S^\alpha$. Autrement dit, le type d'un état de la dynamique $\alpha$ est caractérisé par la relation
\[s\in (typ(s))^\alpha.\]

\paragraph{Dynamiques déterministes.}
La dynamique $\alpha$ est dite \emph{déterministe} (resp. \emph{quasi-déterministe}) si pour toute arête $a\in\overrightarrow{\mathbf{G}}$, la transition $a^\alpha$  est déterministe\footnote{
Voir section \ref{subsec notations transitions} la définition \ref{df transitions quasi deterministes}.
}
(resp. quasi-déterministes).

\subsubsection{Dynamorphismes}\label{subsubsec dynamorphismes}

La définition des morphismes entre dynamiques catégoriques, appelés \emph{dynamorphismes}, telle que nous l'avons écrite dans nos textes d'introduction aux dynamiques catégoriques\footnote{Voir la définition 46 dans \cite{Dugowson:201112} ou la définition 47 dans \cite{Dugowson:201203}.} s'étend immédiatement au cas des dynamiques graphiques.

Ainsi, deux dynamiques $\alpha$ et $\beta$ étant données sur un même graphe $\mathbf{G}$, un \emph{$\mathbf{G}$-dynamorphisme} $\delta:\alpha\looparrowright \beta$ est la donnée,
 pour tout sommet $S$ de $\mathbf{G}$, 
 d'une transition $\delta_S:S^\alpha \rightsquigarrow S^\beta$, autrement dit d'une application $\delta_S:S^\alpha \rightarrow\mathcal{P}(S^\beta)$, telle que 
 pour toute arête $e:S\rightarrow T$ de $\mathbf{G}$, on a
\[\delta_T\odot e^\alpha \subset e^\beta\odot \delta_S, \] où $\odot$ désigne la composée des transitions.

Plus généralement, étant données $\alpha$ une dynamique sur un graphe $\mathbf{F}$ et $\beta$ une dynamique sur un graphe $\mathbf{G}$, un dynamorphisme $(\Delta,\delta):\alpha\looparrowright \beta$ est la donnée
\begin{itemize}
\item d'un morphisme de graphes $\Delta:\mathbf{F}\to \mathbf{G}$,
\item pour tout sommet $S$ de $\mathbf{F}$, 
 d'une transition $\delta_S:S^\alpha \rightsquigarrow (S^\beta)$ telle que 
 pour toute arête $e:S\rightarrow T$ de $\mathbf{F}$, on a
\[\delta_T\odot e^\alpha \subset (\Delta e)^\beta\odot \delta_S, \] où $\odot$ désigne la composée des relations.
\end{itemize}

\subsubsection{Horloge, successions, réalisations}

\begin{df} On appelle \emph{horloge} sur le graphe $\mathbf{G}$ toute dynamique déterministe sur $\mathbf{G}$. Les états d'une horloge sont appelés les \emph{instants} (de cette horloge).
\end{df}

Ainsi, $h$ étant une horloge sur $\mathbf{G}$, tous les ensembles de la forme $a^h(t)$, avec $dom(a)=typ(t)$, sont des singletons. 

\begin{exm}\label{horloge essentielle} \`{A} tout graphe $\mathbf{G}$, on associe canoniquement une horloge, notée $\zeta_\mathbf{G}$ et appelée l'horloge essentielle de $\mathbf{G}$,  en posant
\begin{itemize}
\item pour tout $S\in\dot{\mathbf{G}}$, $\zeta_\mathbf{G}(S)=\{S\}$,
\item pour tout $a\in\overrightarrow{\mathbf{G}}$, $\zeta_\mathbf{a}$ est l'unique application du singleton $\{dom(a)\}$ dans le singleton $\{cod(a)\}$.
\end{itemize}
\end{exm}


\begin{df}\label{df succession}
Nous dirons qu'un instant $t$ de $h$ \emph{succède} à un instant $s$ lorsqu'il existe une arête $a$ telle que $a^h(s)=t$. 
\end{df}

Contrairement à la relation de pré-ordre associée à une horloge définie sur une catégorie\footnote{Voir \cite{Dugowson:201112}, proposition 23.}, la relation binaire de succession définie sur les instants d'une horloge graphique n'est en général ni transitive, ni réflexive (ni bien sûr anti-symétrique, ni symétrique). En particulier, ce n'est pas parce que $t$ succède à $s$ et $s$ à $r$ que $t$ succède à $r$.

\begin{df}\label{df realisation de dyna pour horloge}[Réalisation d'une dynamique pour une horloge]
Une \emph{réalisation} (ou une \emph{solution}) de la dynamique $\alpha:\mathbf{G}\rightarrow Gr(\mathbf{P})$ pour l'horloge $h$ de même moteur $\mathbf{G}$ est un dynamorphisme quasi-déterministe $\sigma : h\looparrowright \alpha$.
\end{df}

\begin{exm} Pour toute dynamique graphique et pour toute horloge de même moteur, on peut toujours définir la \emph{réalisation vide}, qui à tout instant associe l'ensemble vide.
\end{exm}

\subsection{Dynamiques graphiques scandées}

\subsubsection{Définition}

\begin{df}
Une \emph{dynamique graphique scandée sur un graphe $\mathbf{G}$} est un triplet $(\alpha,h,\tau)$ avec 
\begin{itemize}
\item $\alpha$ une dynamique sur $\mathbf{G}$,
\item $h$ une horloge sur  $\mathbf{G}$,
\item $(\tau:\alpha\looparrowright h)$ un $\mathbf{G}$-dynamorphisme déterministe, appelé \emph{scansion} ou \emph{datation} de la dynamique scandée   $(\alpha,h,\tau)$.
\end{itemize}
\end{df}

Puisqu'un dynamorphisme déterministe $\delta:\alpha\looparrowright\beta$  s'identifie à une \textit{application} 
$st(\alpha)\rightarrow st(\beta)$, 
une datation $\tau:\alpha\looparrowright h$ 
associe à tout état $a\in S^\alpha$ un instant $\tau(a)=\tau_S(a)\in S^{h}$, cela pour tout type de temporalité $S\in\dot{\mathbf{G}}$. L'instant $\tau(a)$ sera appelé la \emph{date} de $a$ (pour la datation $\tau$).

On dira parfois d'une telle datation $\tau$ qu'elle est définie \emph{sur} la dynamique $\alpha$. La dynamique scandée $(\alpha,h,\tau)$ sera le plus souvent désignée par sa datation, et l'on parlera ainsi de \emph{la dynamique scandée $\tau:\alpha \looparrowright h$}. S'il n'y a pas d'ambiguïté, on pourra aussi parfois la désigner simplement par $\alpha$.


\begin{exm}\label{exm dynamique essentiellement scandee canonique} \`{A} toute dynamique graphique $\alpha$ sur un graphe $\mathbf{G}$ on associe canoniquement une dynamique scandée par l'horloge essentielle $\zeta_\mathbf{G}$ de $\mathbf{G}$, en munissant $\alpha$ de la datation  $\tau_\alpha:\alpha\looparrowright\zeta_\mathbf{G} $ définie par
\[\forall S\in\dot{\mathbf{G}}, \forall s\in S^\alpha, \tau_\alpha(s)=\{S\}.\]
 On appellera  $\tau_\alpha:\alpha\looparrowright\zeta_\mathbf{G} $ la \emph{dynamique essentiellement scandée canoniquement associée à $\alpha$}.
\end{exm} 

\subsubsection{Dynamorphismes scandés}

Nous donnons la définition suivante sans justification ni exemple, à titre indicatif. En particulier, nous n'expliquerons pas ici la condition de synchronisation entre horloge, d'allure un peu énigmatique --- pourquoi une inclusion dans ce sens et non dans l'autre, et pourquoi pas une égalité ? ---  nous contentant d'indiquer qu'elle se justifie par la notion de réalisation  d'une dynamique scandée lorsqu'une telle réalisation n'est pas définie à tout instant (voir la remarque \ref{rmq inclusion stricte def dyna scand} page \pageref{rmq inclusion stricte def dyna scand}).

\begin{df}[Dynamorphismes scandés]\label{df dynamo scandes} On appelle \emph{dynamorphisme scandé}, ou simplement dynamorphisme, d'une $\mathbf{F}$-dynamique scandée $\rho:\alpha\looparrowright h$ vers une $\mathbf{G}$-dynamique scandée $\tau:\beta\looparrowright k$ la donnée d'un triplet $(\Delta,\delta,d)$ tel que
\begin{enumerate}
\item $(\Delta,\delta)$ est un dynamorphisme de $\alpha$ vers $\beta$,
\item $(\Delta,d)$ est un dynamorphisme de $h$ vers $k$,
\item pour tout $S\in\dot{\mathbf{F}}$, la condition suivante de synchronisation entre $\rho$ et $\tau$ est satisfaite :
\[\tau_{\Delta_S}\odot \delta_S\subset d_S\odot \rho_S. \]  
\end{enumerate}
$\Delta$ sera appelé la partie fonctorielle du dynamorphisme scandé $(\Delta,\delta,d)$, $\delta$ sa partie transitionnelle, et $d$ sa partie horloge.
\end{df}


\mbox{}

En prenant pour flèches les dynamorphismes scandés, on constitue la catégorie des dynamiques graphiques scandées, que nous noterons $\mathbf{DyGraM}$. 

\subsection{Multidynamiques graphiques}

\subsubsection{Définitions}

\begin{df}
Une $\mathbf{G}$-multi-dynamique $\alpha$ consiste en la donnée d'un ensemble non vide $M$, appelé ensemble des \emph{paramètres}\footnote{\label{footnote valeurs parametre} Il aurait sans doute été plus juste de parler de l'ensemble des valeurs du paramètre de $\alpha$, ce paramètre étant l'ensemble $M$ lui-même. Nous ne ferons pas cette nuance.} de la multi-dynamique, et d'un morphisme de graphes $\alpha:\mathbf{G}\longrightarrow Gr(\mathbf{P}^{\underrightarrow{{M}}})$.
\end{df}

Autrement dit, une \emph{multi-dynamique sur le graphe $\mathbf{G}$} peut être vue comme une famille $\alpha=(\alpha_\mu:\mathbf{G}\rightarrow Gr(\mathbf{P}))_{\mu\in M}$ de $\mathbf{G}$-dynamiques, où $M$ est un ensemble non vide, telle que pour tout type de temporalité $T\in\dot{\mathbf{G}}$ il existe un ensemble $T^\alpha$ tel que $\forall \mu\in M, T^{\alpha_\mu}=T^\alpha$.

Pour toute arête $(e:S\rightarrow T)\in\overrightarrow{\mathbf{G}}$ et tout  paramètre $\mu\in M$, on notera indifféremment  $e^{\alpha_\mu}$ ou $e^\alpha_\mu$ la transition de paramètre $\mu$ associée par $\alpha$ à $e$
. 

On peut en particulier voir
$e^\alpha=(e^\alpha_\mu)_{\mu\in M}$
comme une famille indexée par $M$ de transitions $(e^\alpha_\mu:S^\alpha\rightsquigarrow T^\alpha)_{\mu\in M}$. Conformément aux notations introduites en section \ref{subsec notations transitions},  nous désignerons également $e^\alpha$  par la notation
\[e^\alpha: S^\alpha \twosquig_M  T^\alpha,\] ou plus simplement, s'il n'y pas d'ambiguïté sur les paramètres, par
\[e^\alpha: S^\alpha \twosquig  T^\alpha.\]

Les deux expressions  
$(e^\alpha_\mu:S^\alpha\rightsquigarrow T^\alpha)_{\mu\in M}$ et $e^\alpha:S^\alpha \twosquig_M  T^\alpha$ signifierons donc toutes deux la même chose, à savoir que $e^\alpha$ est une famille indexée par $M$ de transitions de $S^\alpha$ vers $T^\alpha$.

Par opposition aux multi-dynamiques, les dynamiques graphiques données par la définition \ref{df dynamiques graphiques} et qui s'identifient aux multi-dynamiques ayant un singleton pour ensemble de paramètres, seront parfois appelées des \textit{mono-dynamiques}.

A noter que l'ensemble $st(\alpha)$ des états d'une multi-dynamique est définie par la même formule que pour une mono-dynamique, à savoir :
\[st(\alpha)=  \bigcup_{S\in\dot{\mathbf{G}}}{S^\alpha}\]

\subsubsection{Multi-dynamorphismes}

La catégorie des multi-dynamiques a pour objets toutes les multi-dynamiques $\alpha:\mathbf{G}\longrightarrow Gr(\mathbf{P}^{\underrightarrow{{M}}})$,
 où $\mathbf{G}$ décrit la classe des graphes et $M$ celle des ensembles, et pour flèches les multi-dynamorphismes définis de la façon suivante.

\begin{df}\label{df multi-dynamorphismes}[Multi-dynamorphismes]
Étant données $\alpha:\mathbf{F}\longrightarrow Gr(\mathbf{P}^{\underrightarrow{{L}}})$ et
$\beta:\mathbf{G}\longrightarrow Gr(\mathbf{P}^{\underrightarrow{M}})$ deux multi-dynamiques, un dynamorphisme $(\theta,\Delta,\delta)$ de $\alpha$ vers $\beta$ consiste en la donnée
\begin{itemize}
\item d'une application $\theta:L\rightarrow M$,
\item d'un morphisme de graphes $\Delta:\mathbf{F}\rightarrow\mathbf{G}$,
\item d'une transition 
$\delta:st(\alpha)\rightsquigarrow st(\beta)$,
\end{itemize} 
\noindent telles que, pour tout $\lambda\in L$, $(\Delta,\delta)$ définit un dynamorphisme de $\alpha_\lambda$ vers $\beta_{\theta(\lambda)}$.
\end{df}

 Autrement dit, pour tout $\lambda\in L$, tous $S$ et $T$ dans $\dot{\mathbf{F}}$ et tout $(e:S\rightarrow T)\in \overrightarrow{\mathbf{F}}$, 
\[\delta_T\odot e^\alpha_\lambda\subset (\Delta e)^\beta_{\theta(\lambda)}\odot\delta_S.\] 

\subsubsection{Quotient paramétrique}

\begin{df}\label{df reduc param multi}
Étant donnée $\alpha:\mathbf{G}\longrightarrow Gr(\mathbf{P}^{\underrightarrow{{M}}})$ une multi-dynamique de moteur $\mathbf{G}$ et d'ensemble de paramètres $M$, et $\sim$ une relation d'équivalence sur $M$, on appelle \emph{quotient de $\alpha$ par $\sim$} et l'on note  $\alpha/{\sim}$ la multi-dynamique $\beta$ sur le même moteur $\mathbf{G}$ et d'ensemble de paramètres $\widetilde{M}=M/{\sim}$ définie par
\begin{itemize}
\item pour tout sommet $S$ de $\mathbf{G}$, on a $S^\beta=S^\alpha$,
\item pour toute arête $(e:S\rightarrow T)\in\overrightarrow{\mathbf{G}}$, pour toute classe $\lambda\in \widetilde{M}$ et tout état $a\in S^\beta$, on a
\[e^\beta_\lambda(a)=\bigcup_{\mu\in\lambda}e^\alpha_\mu(a).\]
\end{itemize}
\end{df}

\subsection{Dynamiques graphiques ouvertes}

\subsubsection{Définition}\label{subsubsec df DyGrOuv}

\begin{df}\label{df DyGrOuv} On appelle \emph{dynamique ouverte (sur $\mathbf{G}$)} toute  \emph{multi-dynamique scandée (sur $\mathbf{G}$)}, autrement dit toute multi-dynamique sur $\mathbf{G}$ munie d'une horloge et d'une datation. Plus précisément, une dynamique ouverte
\[\tau:\alpha=(\alpha_\mu)_{\mu \in M}\looparrowright h\]
sur $\mathbf{G}$, de paramètres $M$ et d'horloge $h$ est constituée
\begin{itemize}
 \item d'une $\mathbf{G}$-multi-dynamique $(\alpha_\mu)_{\mu \in M}$ de paramètres $M$, 
 \item d'une horloge $h$,
 \item  et d'une application $\tau:st(\alpha)\to st(h)$ 
 \end{itemize} 
telles que pour tout $\mu\in M$, $\tau:\alpha_\mu\looparrowright h$ soit une dynamique scandée.
\end{df}

Comme pour les mono-dynamiques scandées, le second membre de l'inclusion 
\[\tau_T \odot e^\alpha_\mu \subset e^h\odot \tau_S\] 
étant toujours un singleton, le seul cas où elle est stricte 
est celui où on l'applique à un état $a\in S^\alpha$ tel que $e^\alpha_\mu(a)=\emptyset$. Par conséquent,
pour tout $\mu\in M$, toute arête $(e:S\rightarrow T)\in\overrightarrow{\mathbf{G}}$, tout état $a\in S^\alpha$ et tout état $b\in e^\alpha_\mu(a)$, on a
\[
\tau_T(b)=e^h (\tau_S (a)).
\]

\subsubsection{Dynamorphismes de dynamiques ouvertes}
\label{subsubsec cat des dynamiques ouvertes}

On constitue la catégorie des dynamiques ouvertes en prenant pour flèches les multi-dynamorphismes scandés ainsi définis :

\begin{df}[Multi-dynamorphismes scandés]\label{df multi-dynamo scandes} On appelle \emph{dynamorphisme ouvert} ou \emph{multi-dynamorphisme scandé}, ou plus simplement  \emph{dynamorphisme}, d'une dynamique ouverte
 $A=({\rho:(\alpha:\mathbf{F}\longrightarrow Gr(\mathbf{P}^{\underrightarrow{L}}))}\looparrowright h)$ 
 vers une autre, 
 $B=({\tau:(\beta:\mathbf{G}\longrightarrow Gr(\mathbf{P}^{\underrightarrow{M}}))}\looparrowright k)$, 
 la donnée d'un quadruplet $(\theta,\Delta,\delta,d)$ tel que
\begin{enumerate}
\item $(\theta,\Delta,\delta)$ est un multi-dynamorphisme de $\alpha$ vers $\beta$,
\item $(\Delta,d)$ est un dynamorphisme de $h$ vers $k$,
\item pour tout $S\in\dot{\mathbf{F}}$, la condition suivante de synchronisation entre $\rho$ et $\tau$ est satisfaite :
\[\tau_{\Delta_S}\odot \delta_S\subset d_S\odot \rho_S. \]  
Étant donné un dynamorphisme $(\theta,\Delta,\delta,d):A\looparrowright B$, on appelle
\begin{itemize}
\item $\theta$ sa  \emph{partie paramétrique},
\item $\Delta$ sa  \emph{partie fonctorielle},
\item $\delta$ sa  \emph{partie transitionnelle},
\item et $d$ sa \emph{partie horloge}.
\end{itemize}
\end{enumerate}
\end{df}

Les mono-dynamiques pouvant être considérées comme des multi-dynamiques particulières et  toute dynamique pouvant être canoniquement scandée\footnote{Voir l'exemple \ref{exm dynamique essentiellement scandee canonique}.}, on vérifie immédiatement que la définition \ref{df multi-dynamo scandes} ci-dessus généralise toutes les définitions de dynamorphismes données précédemment.

\subsubsection{Quotient paramétrique d'une dynamique ouverte}
\label{subsubsec reduc param dyna graph ouverte}

Exactement comme pour les multi-dynamiques (définition \ref{df reduc param multi}), on ainsi définit le quotient paramétrique d'une dynamique graphique ouverte :

\begin{df}\label{df reduc param ouverte}
Étant donnée $\tau:(\alpha:\mathbf{G}\longrightarrow Gr(\mathbf{P}^{\underrightarrow{{M}}}))\looparrowright h)$ une dynamique ouverte de moteur $\mathbf{G}$, d'horloge $h$, de datation $\tau$, d'ensemble de paramètres $M$ et de multi-dynamique $\alpha$, et $\sim$ une relation d'équivalence sur $M$, on appelle \emph{quotient de $\alpha$ par $\sim$} et l'on note  $\alpha/{\sim}$ la dynamique ouverte de même moteur, de même horloge, de même datation, d'ensemble de paramètres $\widetilde{M}=M/{\sim}$  et de multi-dynamique la multi-dynamique $\alpha/{\sim}$ telle qu'elle a été donnée par la définition \ref{df reduc param multi} page \pageref{df reduc param multi}.
\end{df}

\subsection{Réalisations}\label{subsec realisations}

Dans la présente section \ref{subsec realisations}, la notion de réalisation est élargie au cas des multi-dynamiques, des dynamiques scandées et des dynamiques graphiques ouvertes.

\subsubsection{Réalisations d'une dynamique  graphique sur une horloge}

Pour toute $\mathbf{G}$-dynamique $\alpha:\mathbf{G}\rightarrow Gr(\mathbf{P})$ et toute $\mathbf{G}$-horloge $h$, on notera $\mathcal{S}_{(h,\alpha)}$ l'ensemble des $h$-réalisations de $\alpha$, c'est-à-dire, rappelons-le (voir la  définition \ref{df realisation de dyna pour horloge}), l'ensemble des dynamorphismes quasi-déterministes de $h$ dans $\alpha$.

\begin{rmq}\label{rmq convention transitions et fonctions}
En pratique, plutôt que comme un dynamorphisme, nous considérerons en général une $h$-réalisation $\mathtt{a}$ de $\alpha$ comme une fonction $\mathfrak{a}:st(h)\supset df(\mathfrak{a})\dashrightarrow st(\alpha)$, et nous adopterons les conventions suivantes :
\begin{itemize}
\item on ne notera pas différemment la transition $\mathtt{a}$ et la fonction $\mathfrak{a}$, de sorte en particulier que si $t\notin df(\mathfrak{a})$ on s'autorisera à écrire\footnote{Ce qui ne sera pas une source de confusion tant que l'on ne considérera pas l'ensemble vide comme un élément de l'ensemble des états de la dynamique $\alpha$.} $\mathfrak{a}(t)=\emptyset$,
\item aussi, lorsque $t\notin df(\mathfrak{a})$, l'expression $e^\alpha(\mathfrak{a}(t))$ désignera, pour tout $e\in\overrightarrow{\mathbf{G}}$, l'ensemble vide,
\item enfin, pour tout couple $(t_1,t_2)\in st(h)^2$ et tout $e\in\overrightarrow{\mathbf{G}}$, la validité de l'expression
\[
\mathfrak{a}(t_2)\in e^\alpha(\mathfrak{a}(t_1))
\]
sera étendue au cas où $t_2\notin df(\mathfrak{a})$, y compris lorsque $t_1$ n'est pas non plus dans $df(\mathfrak{a})$.
\end{itemize}
\end{rmq}

%
%
%

\subsubsection{Réalisation des dynamiques graphiques scandées}

La notion de \emph{réalisation d'une dynamique scandée} $A=(\tau:\alpha \looparrowright h)$ s'impose très naturellement : il s'agit de toute $h$-réalisation  $\mathtt{a}$ de $\alpha$ telle que $\tau \odot \mathtt{a} \subset 1_{h}$. Autrement dit, pour tout $t\in st(h)$ tel que $\mathtt{a}(t)\neq\emptyset$, on doit avoir $\tau(\mathtt{a}(t))=\{t\}$.

Remarquons que cette notion peut être exprimée directement en termes de dynamorphisme scandé. Pour cela, définissons l'\textit{horloge auto-scandée} associée à une $\mathbf{G}$-horloge $h$ comme étant la dynamique scandée $[h]=( id_h: h \looparrowright h)$, où $id_h$ désigne le dynamorphisme identité de l'horloge $h$, autrement dit
la datation sur $h$ définie pour tout $S\in\dot{\mathbf{G}}$ par ${(id_h)}_S=id_{S^h}$. 

\begin{df}\label{df realisation dynamique scandee} On appelle \emph{réalisation} d'une $\mathbf{G}$-dynamique scandée $A=(\tau:\alpha \looparrowright h)$ tout dynamorphisme scandé quasi-déterministe 
\[(id_\mathbf{G},\mathtt{a}, id_h):[h] \looparrowright (\tau:\alpha \looparrowright h). \] On note $\mathcal{S}_A$ l'ensemble des réalisations de $A$.
\end{df}

Vérifions que la définition \ref{df realisation dynamique scandee} coïncide bien avec la notion attendue, selon laquelle $\mathtt{a}:h\looparrowright \alpha$  est un dynamorphisme quasi-déterministe tel que pour tout $S\in\dot{\mathbf{G}}$, $\tau_S \odot \mathtt{a}_S \subset id_{S^h}$. Or, cette dernière relation est précisément la condition sur les datations que doit satisfaire le triplet $(id_\mathbf{G},\mathtt{a}, id_h)$ pour être un dynamorphisme scandé, les autres conditions étant trivialement satisfaites. Réciproquement, si $(id_\mathbf{G},\mathtt{a}, id_h)$ est un dynamorphisme scandé de  $[h]=( id_h: h \looparrowright h)$ vers $(\tau:\alpha \looparrowright h)$ tel que $\mathtt{a}$ soit quasi-déterministe, la condition sur les datations entraîne que $\mathtt{a}$ est une réalisation de $(\tau:\alpha \looparrowright h)$. 

\begin{rmq}\label{rmq inclusion stricte def dyna scand}
Notons que lorsque la partie transitionnelle $\mathtt{a}$ de la réalisation considérée  n'est pas complète\footnote{C'est-à-dire lorsque la solution n'est pas définie à tout instant.}, l'inclusion exigée sur les datations par la définition \ref{df dynamo scandes} des dynamorphismes scandés peut être stricte, puisqu'il existe dans ce cas un instant 
$t\in S^h$ tel que $\mathtt{a}(t)=\emptyset$, de sorte que
 $\tau_S\odot \mathtt{a}_S (t) 
 =\emptyset \subsetneqq \{t\}
 =id_h\odot id_h (t)$.
\end{rmq}

\begin{rmq}\label{rmq identification des solutions a des fonctions}
Une telle réalisation $(id_\mathbf{G},\mathtt{a}, id_h)$ étant entièrement déterminée par sa partie transitionnelle $\mathtt{a}$, on parlera simplement de la réalisation $\mathtt{a}$. En outre, on identifiera le plus souvent $\mathtt{a}$ avec la \emph{fonction} $\mathfrak{a}=\vert\mathtt{a}\vert$ : une réalisation d'une dynamique scandée sera alors une fonction $\mathfrak{a}$ admettant un domaine de définition $df(\mathfrak{a})\subset st(h)$, et l'on adoptera les conventions de la remarque 
\ref{rmq convention transitions et fonctions} page \pageref{rmq convention transitions et fonctions}.
\end{rmq}

\begin{rmq}
Contrairement à celles des dynamiques autonomes, les réalisations des dynamiques scandées sont, du fait de la relation $\tau(\mathfrak{a}(t))=t$, nécessairement \emph{injectives}.
\end{rmq}

\subsubsection{$h$-réalisations d'une multi-dynamique}
\label{subsubsec h realisation de multi-dynamique}

\begin{df} Étant donnée une $\mathbf{G}$-horloge $h$, une $h$-réalisation  d'une $\mathbf{G}$-multi-dynamique $\alpha:\mathbf{G}\longrightarrow Gr(\mathbf{P}^{\underrightarrow{{L}}})$ est un  $\mathbf{G}$-multi-dynamorphisme quasi-déterministe de $h$ dans $\alpha$.
\end{df}
 
Autrement dit, l'horloge $h$ étant comprise comme multi-dynamique à un seul paramètre, une telle réalisation consiste en le choix d'une valeur du paramètre $\lambda\in L$ et en un $\mathbf{G}$-dynamorphisme quasi-déterministe $\mathfrak{a}$ de $h$ dans $\alpha_\lambda$. Autrement dit, une $h$-réalisation de $\alpha$ est un couple $(\lambda,\mathfrak{a})$. Nous appellerons $\mathfrak{a}$ la \emph{partie externe} de la réalisation $(\lambda,\mathfrak{a})$. Par abus de langage, nous dirons souvent \og la réalisation $\mathfrak{a}$\fg\, au lieu de \og la partie externe $\mathfrak{a}$ d'une réalisation\fg\, : cela ne prêtera jamais à confusion pour la bonne raison que, sauf mention du contraire, nous ne nous intéresserons qu'aux parties externes des réalisations.

Comme d'habitude, plutôt que comme un dynamorphisme, nous considérerons en général une réalisation $\mathfrak{a}$ comme une fonction $\mathfrak{a}:st(h)\supset df(\mathfrak{a})\dashrightarrow st(\alpha)$, et nous utiliserons les conventions indiquées dans la remarque \ref{rmq convention transitions et fonctions} page \pageref{rmq convention transitions et fonctions}.

Notant
\[{\mathcal{S}}_{(h,\alpha)}\]
l'ensemble des $h$-réalisations de la $\mathbf{G}$-multi-dynamique $\alpha=(\alpha_\lambda)_{\lambda\in L}$, on a donc 
\[{\mathcal{S}}_{(h,\alpha)}=
\bigcup_{\lambda\in L} {\mathcal{S}}_{(h,\alpha_\lambda)}.\]

\begin{prop}\label{prop forme des solutions des multi-dynamiques}
Une fonction $\mathfrak{a}:st(h)\dashrightarrow st(\alpha)$ est une $h$-réalisation de la multi-dynamique $\alpha$ si et seulement s'il existe $\lambda\in L$ tel que, pour tout $h$-instant $t\in st(h)$ et tout $e\in \overrightarrow{\mathbf{G}}$ tel que $dom(e)=typ(t)$, on ait
\[ t'\in df(\mathfrak{a}) \Longrightarrow (t\in df(\mathfrak{a}) \,\mathrm{et}\,  \mathfrak{a}(t')\in e^\alpha_\lambda(\mathfrak{a}(t))),\]
où on a posé $t'=e^h(t)$.
\end{prop}

\noindent \textbf{preuve}. Si $\mathfrak{a}$ est une $h$-réalisation de $\alpha$, on a de façon immédiate la condition annoncée par définition d'un dynamorphisme quasi-déterministe et en particulier du fait que si $t\notin df(\mathfrak{a})$, alors $t'\notin df(\mathfrak{a})$, autrement dit $\{\mathfrak{a}(t')\}=\emptyset$, puisque $\{\mathfrak{a}(t')\}\subset e^\alpha_\lambda(\{\mathfrak{a}(t)\})=\emptyset$. 

Réciproquement, supposons qu'il existe $\lambda\in L$ tel que la condition donnée soit satisfaite par la fonction $\mathfrak{a}$. La famille de transitions $(\mathfrak{a}_T=\{\mathfrak{a}_{\vert T}\})_{T\in\dot{\mathbf{G}}}$ --- où $\mathfrak{a}_{\vert T}$ désigne la restriction de la fonction $\mathfrak{a}$ à l'ensemble $T^h$ et où $\{\mathfrak{a}_{\vert T}\}$ désigne la transition associée 
--- vérifie alors, pour tout $e:S\rightarrow T$ et tout $t\in S^h$ :

\begin{itemize}
\item si $t'\notin df(\mathfrak{a})$ et $t\notin df(\mathfrak{a})$, $\mathfrak{a}_T(t')=\emptyset\subset \emptyset=e^\alpha_\lambda(\mathfrak{a}(t))$,
\item si $t'\notin df(\mathfrak{a})$ et $t\in df(\mathfrak{a})$, alors $\mathfrak{a}_T(t')=\emptyset\subset e^\alpha_\lambda(\mathfrak{a}(t))$,
\item si $t'\in df(\mathfrak{a})$, alors $t\in df(\mathfrak{a})$ et $\mathfrak{a}_T(t')=\{\mathfrak{a}(t')\}\subset e^\alpha_\lambda(\mathfrak{a}(t))$, 
\end{itemize}
où $t'$ désigne toujours $e^h(t)$. Ainsi, on a bien, dans tous les cas, $\mathfrak{a}_T \odot e^h \subset e^\alpha_\lambda \odot \mathfrak{a}_S$, de sorte que $\mathfrak{a}$ est un $\mathbf{G}$-dynamorphisme quasi-déterministe de $h$ dans $\alpha_\lambda$, autrement $\mathfrak{a}$ est une $h$-réalisation de $\alpha$.
\begin{flushright}$\square$\end{flushright}

\subsubsection{Réalisation des dynamiques ouvertes}
\label{subsubsec real dyna ouv}

\begin{df} Soit $A=(\tau:(\alpha_\lambda)_{\lambda\in L} \looparrowright h)$ une $\mathbf{G}$-dynamique ouverte, autrement dit une $\mathbf{G}$-multi-dynamique scandée.
On appelle \emph{réalisation de $A$} tout $\mathbf{G}$-multi-dynamorphisme $h$-scandé quasi-déterministe de $[h]=( id_h: h \looparrowright h)$ dans $\tau: \alpha \looparrowright h$. 
\end{df}

\begin{df}Étant donné $\lambda\in L$, on appelle \emph{réalisation de paramètre $\lambda\in L$ de $A$} toute réalisation de la mono-dynamique scandée $(\tau:\alpha_\lambda\looparrowright h)$.
\end{df}

Autrement dit, une  réalisation de $A$ est une réalisation $(\lambda,\mathfrak{a})$ de la multi-dynamique $(\alpha_\lambda)_{\lambda\in L}$ sur l'horloge $h$  qui respecte la datation, \emph{i.e.} telle que
\[\tau(\mathfrak{a}(t))=t.\] 
Conformément aux définitions relatives aux réalisations d'une multi-dynamique graphique, nous dirons que $\mathfrak{a}$ est la partie externe de la réalisation $(\lambda,\mathfrak{a})$.

Ainsi, la partie externe d'une réalisation de $A$ est une réalisation de paramètre $\lambda$ de $A$, pour une certaine valeur du paramètre $\lambda\in L$. Nous noterons   $\mathcal{S}_{(A,\lambda)}$ l'ensemble des réalisations de paramètre $\lambda$ de $A$, et $\mathcal{S}_A$ l'ensemble des parties externes des réalisations de $A$, de sorte que \[\mathcal{S}_A=
\bigcup_{\lambda\in L}{{\mathcal{S}}_{(A,\lambda)}}.\]

En pratique, comme indiqué en section \ref{subsubsec h realisation de multi-dynamique}, nous commettrons souvent l'abus de langage consistant à dire \og la réalisation $\mathfrak{a}$\fg\, de $A$\, au lieu de \og la partie externe $\mathfrak{a}$ d'une réalisation de $A$\fg.

\begin{rmq} Bien entendu, le contenu de la remarque \ref{rmq identification des solutions a des fonctions} page \pageref{rmq identification des solutions a des fonctions} concernant les dynamiques scandées s'applique également aux réalisations des dynamiques ouvertes, et donc en particulier les conventions de la remarque \ref{subsec notations transitions} page \pageref{subsec notations transitions}. De même la notion de \emph{réalisation passant par un état} précisée dans la section \ref{subsubsec realisations passant par etat} s'applique-t-elle également aux réalisations des dynamiques ouvertes.
\end{rmq}

\subsubsection{Réalisations passant par un état}
\label{subsubsec realisations passant par etat}

Nous dirons qu'une réalisation $\mathfrak{a}$ d'une  multi-dynamique graphique $\alpha$ pour une horloge $h$ passe par un état 
$a\in st(\alpha)$ 
s'il existe un instant $t\in st(h)$ tel que 
$\mathfrak{a}(t)=a$. 

Dans le cas des dynamiques scandées, en particulier dans le cas des dynamiques ouvertes, la précision de l'instant $t$ est inutile puisqu'il doit nécessairement être donné par la datation. D'où la définition suivante :

\begin{df}
Étant donnée $A=(\tau:(\alpha_\lambda)_{\lambda\in L}\looparrowright h)$ une dynamique ouverte sur le graphe $\mathbf{G}$, nous dirons qu'une réalisation $\mathfrak{a}$ de $A$ passe par un état $a\in st(\alpha)$, et nous écrirons
\[\mathfrak{a}\rhd a,\] si et seulement si $\mathfrak{a}(\tau(a))=a$.
\end{df}

Dans le cas où $a$ et $b$ sont deux états de la dynamique ouverte $A$ tels que $\tau(b)$ succède\footnote{Voir la définition \ref{df succession}.} à $\tau(a)$, nous écrirons
\[\mathfrak{a}\rhd a, b\] pour exprimer que $\mathfrak{a}$ passe par $a$ \emph{puis} qu'elle passe par $b$.

\section{Interactions et familles dynamiques}
\label{section familles dynamiques}

Nous définissons ici les familles dynamiques, c'est-à-dire les familles de dynamiques graphiques ouvertes en interaction grâce à une relation liant les membres de la famille, membres parmi lesquels une sorte de chef d'orchestre est chargé de donné une temporalité commune de référence de façon à ce qu'une nouvelle dynamique graphique ouverte soit produite par la famille qui se trouvera ainsi à son tour en situation de nouer des relations productives avec d'autres dynamiques ouvertes.

Comme nous avons commencé à le faire dans les sections précédentes, nous utiliserons les lettres $A$, $B$, etc., pour désigner des dynamiques ouvertes :
\[A=({\tau:(\alpha:\mathbf{G}\longrightarrow Gr(\mathbf{P}^{\underrightarrow{L}}}))\looparrowright h),\]
\[B=({\rho:(\beta:\mathbf{F}\longrightarrow Gr(\mathbf{P}^{\underrightarrow{M}}))}\looparrowright k),\quad \mathrm{etc.}\]

Ici, $A$ est par exemple une $\mathbf{G}$-dynamique ouverte, d'ensemble de paramètres $L$, définie par la famille de $\mathbf{G}$-dynamiques $(\alpha_\lambda:\mathbf{G}\rightarrow Gr(\mathbf{P}))_{\lambda\in L}$, d'horloge $h$ et scandée par le $\mathbf{G}$-dynamorphisme déterministe $\tau:\alpha\looparrowright h$. 
S'il n'y a pas d'ambiguïté sur le graphe $\mathbf{G}$, on pourra écrire de façon plus concise
\[A=(\tau:(\alpha_\lambda)_{\lambda\in L} \looparrowright h).\]

De même, une famille indexée par un ensemble $I$ de dynamiques graphiques ouvertes pourra-t-elle être désignée par $(A_i)_{i\in I}$ avec
\[A_i=({\tau_i:(\alpha_i:\mathbf{G}_i\longrightarrow Gr(\mathbf{P}^{\underrightarrow{L_i}}}))\looparrowright h_i).\]

Conformément aux notations de la section \ref{subsubsec real dyna ouv}, $\mathcal{S}_{A_i}$ désignera alors, pour tout $i\in I$, l'ensemble des parties externes des réalisations de la dynamique ouverte $A_i$.

Dans les définitions suivantes, on suppose que les $A_i$ suivent les notations ci-dessus.

\begin{df}\label{df interaction} [Interaction] Étant donnés $I$ un ensemble non vide et  $(A_i)_{i\in I}$  une famille indexée par $I$ de dynamiques ouvertes avec 
\[A_i=({\tau_i:(\alpha_i:\mathbf{G}_i\longrightarrow Gr(\mathbf{P}^{\underrightarrow{L_i}}}))\looparrowright h_i),\] on appelle \emph{relation binaire multiple \og réalisations/paramètres\fg\, pour cette famille} toute relation binaire multiple $R\in \mathcal{BM}_{(\mathcal{S},\mathcal{L})}$ non vide avec en entrée $\mathcal{S}$ et en sortie $\mathcal{L}$ respectivement données par
\[\mathcal{S}=(\mathcal{S}_{A_i})_{i\in I} 
\quad\mathrm{et}\quad
\mathcal{L}=(L_i)_{i\in I}.
\] En outre, une telle relation binaire multiple \og réalisations/paramètres\fg\, $R$ pour la famille $(A_i)_{i\in I}$ est dite \emph{cohérente} si pour tout $(\mathfrak{a}_i,\lambda_i)_{i\in I}\in R$, on a 
\[\mathfrak{a}_i\in \mathcal{S}_{(A_i,\lambda_i)}.\]
Une relation binaire multiple \og réalisations/paramètres\fg\, cohérente pour une famille de dynamiques graphiques ouvertes sera également appelée une \emph{interaction} pour cette famille.
\end{df}

\pagebreak[3]

\begin{df}\label{df famille dynamique}[Familles dynamiques] On appelle \emph{famille dynamique} (au sens graphique) la donnée $(I,i_0, (A_i)_{i\in I}, R, (\Delta_i,\delta_i)_{i\neq i_0})$
\begin{itemize}
\item d'un ensemble $I$ non vide, appelé \emph{index} de la famille,
\item d'un élément $i_0\in I$, appelé \emph{indice synchronisateur} de la famille,
\item d'une famille indexée par $I$ de dynamiques ouvertes $(A_i)_{i\in I}$  appelées \emph{composantes} de la famille dynamique, 
\item d'une interaction\footnote{Voir ci-dessus la définition \ref{df interaction}.} $R\in \mathcal{BM}_{(\mathcal{S},\mathcal{L})}$ pour la famille $(A_i)_{i\in I}$, 
\item d'une famille 
\[
((\Delta_i,\delta_i) 
: h_{i_0}
\looparrowright 
h_i)_{i\in I\setminus\{i_0\}}
\] 
de dynamorphismes déterministes, appelés \emph{synchronisations}.
\end{itemize}
\end{df}

\begin{df} La \emph{structure connective d'une famille dynamique} est la structure connective de sa relation binaire multiple. L'\emph{ordre
connectif d'une famille dynamique} est l'ordre connectif\footnote{Sur la notion d'ordre connectif, voir par exemple, dans le cas fini, la définition 16 de \cite{Dugowson:201012}, et dans le cas général, la section § 1.11 de \cite{Dugowson:201112} et \cite{Dugowson:201203}.} de sa structure connective.
\end{df}

\section{Dynamiques ouvertes engendrées par une famille dynamique}
\label{section DyGrOuv engendrees}

\`{A} toute famille dynamique graphique $\mathcal{F}$, nous allons associer plusieurs dynamiques graphiques ouvertes que l'on considérera comme produites par cette famille, les différences entre elles portant uniquement sur leurs paramétrisations. Définie en section \ref{subsec dyna primo}, la première de ces dynamiques, appelée la dynamique primo-engendrée par la famille dynamique considérée, est celle pour laquelle la paramétrisation la plus large possible est utilisée. Les trois autres dynamiques que nous considérerons --- la dynamique fonctionnellement engendrée (section \ref{subsubsec dyna fonct eng}), la dynamique souplement engendrée (section \ref{subsubsec dyna soupl eng})  et la  dynamique mono-engendrée (section \ref{subsubsec dyna mono eng}) --- s'exprimeront comme quotients paramétriques de la dynamique primo-engendrée.

Rappelons que, conformément aux remarques \ref{rmq notation binaire rbC a} et \ref{rmq transitions et relations binaires}, l'image $\mathrm{Im}(rb(R))$ de la relation binaire $rb(R)$ associée à une relation binaire multiple $R\in \mathcal{BM}_{(\mathcal{S}, \mathcal{L})}$ avec  $\mathcal{S}=(\mathcal{S}_{A_i})_{i\in I}$ et $\mathcal{L}= (L_i)_{i\in I}$  est définie par 
\[\mathrm{Im}(rb(R))
=
\{
(\lambda_i)_{i\in I}\in \Pi_I(L_i),
\exists (\mathfrak{a}_i)_{i\in I}\in \Pi_I(\mathcal{S}_{A_i}),
(\mathfrak{a}_i,\lambda_i)_{i\in I}\in R
\}
.\]
Rappelons également\footnote{Voir la remarque \ref{rmq notation rbBmoins1} page \pageref{rmq notation rbBmoins1}.} que pour toute relation binaire $B:\mathcal{A}\rightsquigarrow \Lambda$, et pour tout $\mu\in \Lambda$, nous notons 
\[B^{-1}(\mu)=\{\mathfrak{a}\in\mathcal{A},B(\mathfrak{a},\mu)\}.\] En particulier, pour toute relation binaire multiple $R\in \mathcal{BM}_{(\mathcal{S}, \mathcal{L})}$ avec $\mathcal{S}=(\mathcal{S}_{A_i})_{i\in I}$ et $\mathcal{L}= (L_i)_{i\in I}$, et pour tout $\mu\in\Pi_I(L_i)$, l'expression $rb(R)^{-1}(\mu)$ désigne l'ensemble
\[rb(R)^{-1}(\mu)=\{\mathfrak{a}\in\Pi_I(\mathcal{S}_{A_i}), (\mathfrak{a},\mu)\in rb(R)\}.\]

Rappelons enfin qu'une expression de la forme 
\[\mathfrak{a}\triangleright a, b\] se lit \og $\mathfrak{a}$ passe par $a$ puis par $b$\fg, et que la signification en est donnée en section \ref{subsubsec realisations passant par etat}. Ces rappels faits, nous pouvons maintenant poser la définition de la dynamique ouverte produite par une famille dynamique.

\subsection{Dynamique $[\mathcal{F}]_\mathrm{p}$, primo-engendrée par $\mathcal{F}$}
\label{subsec dyna primo}

\begin{df}\label{df dynamique primo-engendree} Étant donnée  $\mathcal{F}=(I,0, (A_i)_{i\in I}, R, (\Delta_i,\delta_i)_{i\neq i_0})$ une famille dynamique d'indice synchronisateur noté $0$ et ayant pour composantes les dynamiques ouvertes
$A_i=({\tau_i:(\alpha_i:\mathbf{G}_i\longrightarrow Gr(\mathbf{P}^{\underrightarrow{L_i}}}))\looparrowright h_i)$,
on appelle \emph{dynamique primo-engendrée par $\mathcal{F}$}, la dynamique ouverte notée $[\mathcal{F}]_\mathrm{p}$ définie par  
\[[\mathcal{F}]_\mathrm{p}=({\rho:(\beta:\mathbf{F}\longrightarrow Gr(\mathbf{P}^{\underrightarrow{M}}))}\looparrowright k)\]
avec :
$\mathbf{F}=\mathbf{G}_0$,
$k=h_0$,
$M=\mathrm{Im}(rb(R))$,
$\beta$ est la multi-dynamique graphique sur $\mathbf{F}$ d'ensemble de paramètres $M$ définie
pour tout sommet $S\in \dot{\mathbf{F}}$ par\footnote{\label{footnote ordre des termes} Dans le produit cartésien écrit ci-dessous, $S^{\alpha_0}\times \prod_{i\neq 0}(\Delta_i S)^{\alpha_i}$, l'ordre des termes ne joue aucun rôle (l'ensemble $I$ n'étant pas \emph{a priori} lui-même ordonné) et nous aurions pu aussi bien écrire $\prod_{i\neq 0}(\Delta_i S)^{\alpha_i}\times S^{\alpha_0}$ pour désigner le même ensemble, à savoir l'ensemble des familles d'éléments $(a_i)_{i\in I}$ indexées par $I$ et prenant leurs valeurs dans les ensembles indiqués.}
\[S^\beta
=
\{(a_i)_{i\in I}\in S^{\alpha_0}\times \prod_{i\neq 0}(\Delta_i S)^{\alpha_i},
\forall i\neq 0, \tau_{i(\Delta_i S)}(a_i)=\delta_i(\tau_{0(S)}(a_0))\},\]
et pour toute arête $(e:S\rightarrow T)\in\overrightarrow{\mathbf{F}}$, tout paramètre $\mu\in M$ et tout état  $a=(a_i)_{i\in I}\in S^\beta$ par 
\[ 
e^\beta_\mu(a)=
\{
b\in T^\beta, 
\tau_{0(T)}(b_0)=e^{h_0}(\tau_{0(S)}(a_0))\,
\mathrm{et}\,
\exists (\mathfrak{a}_i)_{i\in I}\in  rb(R)^{-1}(\mu),
\forall i\in I, \mathfrak{a}_i\triangleright a_i, b_i
\}
\]
et, enfin, $\rho$ est la datation définie pour tout $a=(a_i)_{i\in I}\in S^\beta$ par $\rho_{(S)}(a)=\tau_{0(S)}(a_0)$.
\end{df}


\subsection{Trois quotients paramétriques de $[\mathcal{F}]_\mathrm{p}$}\label{subsec trois reduc pour F}

Comme indiqué au début de cette section \ref{section DyGrOuv engendrees}, nous allons associer trois autres dynamiques graphiques à la famille dynamique $\mathcal{F}$, obtenues comme quotients paramétrique de $[\mathcal{F}]_\mathrm{p}$. La raison en est que l'ensemble $M$ des paramètres de $[\mathcal{F}]_\mathrm{p}$ est en général \og trop gros\fg\, en ce sens que bien souvent le choix d'une valeur quelconque dans $M$ ne sera pas compatible avec le libre fonctionnement de la dynamique engendrée et, de ce fait, pourrait sembler peu naturel.

Chacune des trois dynamiques définies ci-après --- la dynamique fonctionnellement engendrée $[\mathcal{F}]_\mathrm{f}$ (section \ref{subsubsec dyna fonct eng}), la dynamique souplement engendrée $[\mathcal{F}]_\mathrm{s}$ (section \ref{subsubsec dyna soupl eng})  et la  dynamique mono-engendrée $[\mathcal{F}]_\mathrm{m}$ (section \ref{subsubsec dyna mono eng}) ---  sont des quotients paramétriques\footnote{Voir la section \ref{subsubsec reduc param dyna graph ouverte}.} de la dynamique primo-engendrée $[\mathcal{F}]_\mathrm{p}$ par une certaine relation d'équivalence sur l'ensemble $M$ des paramètres de celle-ci, le principe de construction de cette relation d'équivalence étant le même dans les trois cas, à savoir qu'elle résulte du choix de ce que nous appellerons une famille de tas paramétriques.

\subsubsection{Tas paramétriques et équivalence sur $M$}\label{subsubsec tas param et equiv sur M}
Pour chaque $i\in I$, ayant fixé une certaine partie $N_i\subset L_i$  appelée le \emph{tas d'indice $i$}, partie intuitivement destinée à rassembler les valeurs du paramètre\footnote{Voir la note \ref{footnote valeurs parametre}.} de la dynamique $A_i$ dont on considère qu'il appartient au libre fonctionnement de la dynamique engendrée que de les déterminer, et ayant ainsi constitué une famille $\mathcal{N}=(N_i)_{i\in I}$ de tas paramétriques, on considère sur $M\subset \Pi_I\mathcal{L}$ la relation d'équivalence $\sim_\mathcal{N}$ définie, pour tout couple $((\lambda_i)_{i\in I},(\lambda'_i)_{i\in I})\in M^2$, par
\[
(\lambda_i)_{i\in I}\sim_\mathcal{N}(\lambda'_i)_{i\in I}
\]
\[\Leftrightarrow\]
\[
\forall i\in I, \,
(\lambda_i=\lambda'_i)
\,\mathrm{ou}\,
(\lambda_i\in N_i\ni \lambda'_i).\]

La dynamique engendrée par $\mathcal{F}$ au sens des tas paramétriques $\mathcal{N}$ est alors\footnote{Voir la section \ref{subsubsec reduc param dyna graph ouverte}.}
\[[\mathcal{F}]_\mathcal{N}=[\mathcal{F}]_\mathrm{p}/\sim_\mathcal{N}.\]

\subsubsection{Familles $R$-compatibles} Pour préciser la manière dont les tas paramétriques sont définis dans les constructions des sections suivantes, nous aurons besoin de faire appel à la notion de famille compatible pour la relation binaire multiple $R$. \`{A} une telle relation $R\in\mathcal{BM}_{(\mathcal{S},\mathcal{L})}$  nous avons associé en section \ref{subsubsec trois injections} une relation    multiple $rd(R)$ d'index $2I=I\sqcup I=I\times\{0,1\}$. Pour toute partie $W\subset 2I$ et toute famille $(r_w)_{w\in W}$ d'éléments pris respectivement --- selon que $w$ est de la forme $(i,0)$ ou $(i,1)$ ---  dans les ensembles constituant la famille $\mathcal{S}$ ou ceux de la famille $\mathcal{L}$, nous dirons que cette famille $(r_w)_{w\in W}$  est  \emph{$R$-compatible} (ou \emph{compatible} avec $R$) si elle est la restriction à $W$ d'une famille appartenant au graphe de $rd(R)$. Autrement dit, prenant $\mathcal{S}=(\mathcal{S}_{A_i})_{i\in I}$, $\mathcal{L}=(L_i)_{i\in I}$, posant $\mathcal{E}=(\mathcal{S},\mathcal{L})$, notant
\footnote{\label{footnote ordre des facteurs indifferents}Sans tenir compte de l'ordre des facteurs.}
\[
\Pi_W\mathcal{E}=\prod_{(i,0)\in W}{\mathcal{S}_{A_i}}\times \prod_{(i,1)\in W}{L_i}
\] et appelant $\Lambda_W$ la restriction 
\[\Lambda_W:\Pi_{2I}\mathcal{E} \rightarrow \Pi_{W}\mathcal{E} \] définie par
\[\Lambda_W((r_j)_{j\in 2I})=(r_w)_{w\in W},\]
 on a
\[(r_w)_{w\in W}\in \Pi_{W}\mathcal{E}\,\mathrm{est}\, R\mathrm{-compatible} \]
\[\Leftrightarrow\]
\[(r_w)_{w\in W}\in \Lambda_W(rd(R)),\] où $rd(R)\subset \Pi_{2I}\mathcal{E}$ désigne le graphe de la relation multiple $rd(R)$.

Par exemple, un paramètre $l_k\in L_k$ sera dit $R$-compatible s'il existe $\mu\in M$, où $M$ est l'ensemble de paramètres donné par la définition \ref{df dynamique primo-engendree} ci-dessus, tel que $l_k=\mu_k$.

De plus, étant donnés $X\subset 2I$, $Y\subset 2I$  et deux familles $q=(q_x)_{x\in X}\in \Pi_X\mathcal{E}$ et $r=(r_y)_{y\in Y}\in \Pi_Y\mathcal{E}$  compatibles entre elles au sens où pour tout $z\in X\cap Y$ on a $q_z=r_z$ --- ce qui est le cas notamment si $X\cap Y=\emptyset$ --- nous noterons $q+r$ la famille $q+r=s=(s_z)_{z\in X\cup Y}$ telle que, pour tout $z\in X\cup Y$, on a $z\in X \Rightarrow s_z=q_z$ et $z\in Y \Rightarrow s_z=r_z$. Bien entendu\footnote{Conformément aux notes \ref{footnote ordre des termes} et \ref{footnote ordre des facteurs indifferents}. Voir aussi la proposition 1 de l'article \cite{Dugowson:201505}.}, $q+r=r+q$.  Nous utiliserons en particulier cette notation avec des familles de la forme \[(\mathfrak{a}_j)_{j\in J\subset I}\in \Pi_J(\mathcal{S}),\] sans les ré-écrire comme il le faudrait en toute rigueur sous la forme \[(\mathfrak{a}_j)_{(j,0)\in J\times \{0\}\subset 2I}\in \Pi_{J\times \{0\}}\mathcal{E},\] de même qu'avec des familles de la forme \[(l_k)_{k\in K\subset I}\in \Pi_K(\mathcal{L})\] comprises implicitement comme désignant \[(l_k)_{(k,1)\in K\times \{1\}\subset 2I}\in \Pi_{K\times \{1\}}\mathcal{E}.\] Par exemple, dire que l'on a $l=(l_i)_{i\in I}\in M$, ce qui revient à dire que la famille $l=(l_i)_{i\in I}$ est $R$-compatible, équivaut encore à
\[\exists \mathfrak{a}\in \Pi_I(\mathcal{S}), l+\mathfrak{a}\in rd(R). \]

\subsubsection{Dynamique $[\mathcal{F}]_\mathrm{f}$ fonctionnellement engendrée par $\mathcal{F}$}\label{subsubsec dyna fonct eng}

Pour tout $k\in I$, on définit le tas fonctionnel $N^{\mathrm{f}}_k\subset L_k$ de la façon suivante : un élément $l_k\in L_k$ vérfie $l_k\in N^{\mathrm{f}}_k$ si et seulement si $l_k$ est $R$-compatible\footnote{\label{footnote condition l compatible}Logiquement, la condition \og $l_k$ est $R$-compatible\fg\, ne change rien, puisqu'au bout du compte la relation d'équivalence définie par les tas le sera sur l'ensemble $M$, mais c'est sans doute plus clair de se limiter aux $l_k$ effectivement concernés.} et si l'implication suivante est satisfaite
\[
\forall \mathfrak{a}\in\Pi_{i\neq k}{\mathcal{S}_{A_i}},
\forall \lambda_k\in L_k,
\left.
\begin{tabular}{c}
$\mathfrak{a}+ l_k$ est $R$-compatible\\ 
$\mathfrak{a}+ \lambda_k$ est $R$-compatible\\ 
\end{tabular} 
\right\rbrace
\Rightarrow
l_k=\lambda_k.
\]

Intuitivement, on place dans le tas $N^{\mathrm{f}}_k$ les paramètres de la dynamique graphique ouverte $A_k$ dont la valeur est déterminée \emph{via} l'interaction $R$ par les réalisations des \emph{autres} dynamiques en jeu.

La famille $\mathcal{N}^{\mathrm{f}}=(N^{\mathrm{f}}_k)_{k\in I}$ des tas paramétriques fonctionnels étant ainsi définie, on applique le procédé décrit dans la section \ref{subsubsec tas param et equiv sur M}, obtenant ainsi une relation d'équivalence $\sim_{\mathrm{f}}$ sur $M$, et on appelle \emph{dynamique ouverte graphique fonctionnellement engendrée par la famille dynamique graphique $\mathcal{F}$}, et l'on note $[\mathcal{F}]_\mathrm{f}$, la dynamique graphique ouverte
\[ [\mathcal{F}]_\mathrm{f}=[\mathcal{F}]_\mathrm{p}/\sim_{\mathrm{f}},\] d'ensemble de paramètres $M/\sim_{\mathrm{f}}$.

\subsubsection{Dynamique $[\mathcal{F}]_\mathrm{s}$ souplement engendrée par $\mathcal{F}$}\label{subsubsec dyna soupl eng}

Pour tout $k\in I$, on définit le tas fonctionnel $N^{\mathrm{s}}_k\subset L_k$ comme l'ensemble des éléments bloqués de $L_k$, un élément de $L_k$ étant dit \emph{bloqué} (pour l'interaction $R$) s'il n'est pas libre, tandis qu'un élément $\lambda_k\in L_k$ est dit \emph{libre} ou \emph{souple} s'il est $R$-compatible\footnote{La remarque de la note \ref{footnote condition l compatible} ci-dessus s'applique encore ici, et s'appliquerait aussi bien aux éléments bloqués.} et si quel que soit $\mathfrak{a}_k\in\mathcal{S}_{A_k}$, quel que soit $\mu\in \Pi_{j\neq k}L_j$ et quel que soit $\mathfrak{b}\in \Pi_{j\neq k}\mathcal{S}_{A_j}$ on a l'implication
\[
\left.
\begin{tabular}{c}
$\lambda_k+\mathfrak{a}_k+\mu$ est $R$-compatible\\ 
$\mu+\mathfrak{b}$ est $R$-compatible\\ 
\end{tabular} 
\right\rbrace
\Rightarrow
\lambda_k+\mu+\mathfrak{a}_k+\mathfrak{b}\,\mathrm{est}\, R\mathrm{-compatible}.
\]

La famille $\mathcal{N}^{\mathrm{s}}=(N^{\mathrm{s}}_k)_{k\in I}$ constituée des tas de paramètres bloqués étant ainsi définie, on applique le procédé décrit dans la section \ref{subsubsec tas param et equiv sur M}, obtenant ainsi une relation d'équivalence $\sim_{\mathrm{s}}$ sur $M$, et on appelle \emph{dynamique ouverte graphique souplement engendrée par la famille dynamique graphique $\mathcal{F}$}, et l'on note $[\mathcal{F}]_\mathrm{s}$, la dynamique graphique ouverte
\[ [\mathcal{F}]_\mathrm{s}=[\mathcal{F}]_\mathrm{p}/\sim_{\mathrm{s}},\] d'ensemble de paramètres $M/\sim_{\mathrm{s}}$.

Intuitivement, sont mis dans les tas de paramètres bloqués ceux dont le choix par un agent extérieur à la famille dynamique considérée pourrait être en retour remis en cause par ce que nous pourrions appeler le libre fonctionnement de cette famille, représenté ici par $\mathfrak{a}_k$, $\mu\in \Pi_{j\neq k}L_j$ et  $\mathfrak{b}$. Cette formalisation s'éclairera, nous l'espérons, sur des exemples\footnote{En l'occurrence, nous aurons notamment été guidé par l'examen de ce qui doit être considéré comme libre paramètre ou non dans le cas d'un ressort soumis à différents jeux de contraintes tels que la connaissance du comportement du ressort permette de savoir à quel jeu de contrainte il est soumis...} qui seront donnés dans des articles ultérieurs. 

\subsubsection{Dynamique $[\mathcal{F}]_\mathrm{m}$ mono-engendrée par $\mathcal{F}$}\label{subsubsec dyna mono eng}

On prend pour relation d'équivalence sur $M$ la relation d'équivalence maximale $\sim_{\mathrm{m}}$, de sorte que $M/\sim_{\mathrm{m}}$ est réduit à un point, et on appelle \emph{mono-dynamique graphique  engendrée par la famille dynamique graphique $\mathcal{F}$}, ou encore \emph{dynamique graphique  mono-engendrée par la famille dynamique graphique $\mathcal{F}$}, et l'on note $[\mathcal{F}]_\mathrm{m}$, la dynamique graphique scandée
\[ [\mathcal{F}]_\mathrm{m}=[\mathcal{F}]_\mathrm{p}/\sim_{\mathrm{m}},\] qui ne dépend d'aucun paramètre. Bien entendu, cette mono-dynamique scandée peut toujours être vue comme une dynamique ouverte, l'ensemble des valeurs prises par le paramètre\footnote{Sur notre usage du mot paramètre, voir la note \ref{footnote valeurs parametre}.} se réduisant à un singleton. Cette construction revient à prendre pour tas paramétriques d'indice $k$ l'ensemble $L_k$ lui-même ou, de façon équivalente\footnote{Comme indiqué dans la note \ref{footnote condition l compatible}.},  l'ensemble des valeurs de $L_k$ qui sont $R$-compatibles.

\section{Et maintenant ?}\label{section conclusion}

Comme annoncé en introduction, seules les définitions ont été \og parachutées\fg\, dans le présent article, et pour en illustrer la signification de nombreux exemples de dynamiques graphiques ouvertes devront être présentés ultérieurement. Par ailleurs, l'étude de la signification et des relations logiques entre les différents types d'engendrement dynamique présentés en section \ref{section DyGrOuv engendrees} reste à ce stade très largement à explorer.
Quoi qu'il en soit,
les dynamiques considérées ici auraient pu aussi bien être qualifiées de \og pré-dynamiques\fg, dans la mesure où la restriction à des moteurs graphiques signifie que l'aspect proprement dynamique des choses, fondé sur l'enchaînement des écoulements temporels, n'est pas encore envisagé à ce stade. Métaphoriquement, nous pourrions dire que si les dynamiques sont vues comme les ricochets d'une pierre plate lancée sur un plan d'eau, les dynamiques graphiques ne s'intéressent qu'à un seul rebond : les choses réellement intéressantes ne commencent qu'avec le désir de produire le meilleur enchaînement possible de rebonds...  Aussi, les notions ici définies sont-elles en fait destinées à servir de socle aux notions de \emph{dynamique catégorique ouverte}, de \emph{famille dynamique} (au sens des dynamiques catégoriques) et de \emph{dynamique ouverte produite par une telle famille}. Ce socle aura été rendu nécessaire par le type de problèmes de recollement évoqué dans la conférence \cite{Dugowson:20150506}, qui conduit en effet à ce qu'en général la dynamique ouverte produite par une famille de dynamiques \emph{catégoriques} en interaction ne soit pas elle-même catégorique. Par contre, sur la base des définitions posées dans la présente note, nous sommes assurés que la dynamique ainsi produite est au moins une dynamique graphique ouverte. Nous verrons dans l'article suivant\footnote{\og Dynamiques sous-catégoriques ouvertes en interaction (définitions et théorème de stabilité)\fg, \cite{Dugowson:20150809}.} qu'il est  possible d'élargir le socle des dynamiques graphiques à ce que nous appellerons les \emph{dynamiques sous-catégoriques}, une famille de dynamiques sous-catégoriques ouvertes en interaction produisant (au sens de l'un des types d'engendrement déjà considéré ci-dessus) une dynamique ouverte encore sous-catégorique. Par conséquent si, comme annoncé dans l'introduction du présent article, les \emph{dynamiques catégoriques ouvertes}, qui y seront définies au passage, engendrent par leurs interactions des dynamiques qui ne sont plus nécessairement catégoriques, du moins seront-elles sous-catégoriques.

\bibliographystyle{plain}



\tableofcontents

\end{document}